\documentclass[12pt,twoside]{article}
\usepackage{amsmath,amsbsy,amsfonts}
\def\nd{\noindent}
\def\thend{\rule{3mm}{3mm}}

\newtheorem{theorem}{Theorem}[section]
\newtheorem{definition}{Definition}[section]
\newtheorem{prop}{Proposition}[section]

\newtheorem{remark}{Remark}[section]

\begin{document}

\setlength{\baselineskip}{6.5mm} \setlength{\oddsidemargin}{8mm}
\setlength{\topmargin}{-3mm}
\title{\Large\sf On local compactness in quasilinear elliptic problems}
\author{\sf K. Adriouch \,\, and \,\, A. El Hamidi  }
\date{}
\maketitle

{\bf \sf Keywords.} local compactness, quasilinear systems, Nehari manifold, critical level 

{\bf \sf Mathematics subject Classification.} 35J55; 35J60

\begin{abstract}
One of the major difficulties in nonlinear elliptic problems involving critical nonlinearities is the compactness of 
Palais-Smale sequences. In their celebrated work \cite{BN}, Br\'ezis and Nirenberg introduced the notion of 
critical level for these sequences in the case of a critical perturbation of the Laplacian homogeneous eigenvalue problem. 
In this paper, we give a natural and general formula of the critical level for a large class of nonlinear elliptic critical problems. 
The sharpness of our formula is established by the construction of suitable Palais-Smale sequences which are not relatively compact. 
\end{abstract}
\section{Introduction}
In nonlinear elliptic variational problems involving critical nonlinearities, one of the major difficulties 
is to recover the compactness of Palais-Smale sequences of the associated Euler-Lagrange functional. 
Such questions were first studied, in our knowledge, by Br\'ezis and Nirenberg in their well-known work \cite{BN}. 
The concentration-compactness principle due to Lions \cite{PLL} is widely used 
to overcome these difficulties. Other methods, based on the convergence almost everywhere of the gradients of Palais-Smale sequences, 
can be also used to recover the compactness. We refer the reader to the papers by Boccardo and Murat \cite{BM} and by  
J. M. Rakotoson \cite{RAK2} for bounded domains. For arbitray domains, we refer to the recent work by A. El Hamidi and 
J. M. Rakotoson \cite{ER}.

In \cite{BN}, the authors studied the critical perturbation of the eigenvalue problem: 
\begin{equation} \label{bn_eq}
\left\{
\begin{array}{ccc}
-\Delta u &=& \lambda u + u^{2^*-1}  \;\; {\rm in } \;\; \Omega , \\ 
u &>& 0 \;\; {\rm in } \;\; \Omega , \\ 
u &=& 0 \;\; {\rm on } \;\; \partial \Omega,
\end{array}
\right.
\end{equation}
where $\Omega$ is a bounded domain in $\mathbb{R}^N$, $N \geq 3$, with smooth boundary, $2^*=\frac{2N}{N-2}$ 
is the Sobolev critical exponent of the embedding $W^{1,2}(\Omega) \subset L^p(\Omega)$, 
and $\lambda$ is a positive parameter. The authors introduced an important condition on the level corresponding 
to the energy of Palais-Smale sequences which guarantees their relative compactness. Indeed, let $(u_n)$ be a 
Palais-Smale sequence for the Euler-Lagrange functional  
$$
I_{\lambda} (u) = \frac{1}{2} \int_{\Omega} |\nabla u|^2  - \frac{\lambda}{2} \int_{\Omega} |u|^2 - 
\frac{1}{2^*} \int_{\Omega} |u|^{2^*}.
$$
More precisely, the authors showed that if 
\begin{equation} \label{seuil}
\lim_{n \to +\infty} I_{\lambda} (u_n) < \frac{1}{N} S^{\frac{N}{2}} 
\end{equation}
then $(u_n)$ est relatively compact, which implies the existence of nontrivial critical points of $I_\lambda$. Here, $S$ 
denotes the best Sobolev constant in the embedding $W^{1,2}_0(\Omega) \subset L^{2^*}(\Omega)$. In this work, 
we begin by giving the generalization of condition (\ref{seuil}) for the quasilinear equation 
\begin{equation}
\begin{array}{l}
-\Delta_p u = \lambda f(x,u) + |u|^{p^*-2}u  \;\; {\rm in } \;\; \Omega , \\ \\
u|_{\Gamma} = 0 \;\; \mbox{and} \;\; \frac{\partial{u}}{\partial \nu}|_{\Sigma}=0,
\end{array}
\end{equation}
where $\Omega$ is a bounded domain in $\mathbb{R}^N$, $N \geq 3$, 
with smooth boundary $\partial \Omega = \overline \Gamma \cup \overline \Sigma$,  
where $\Gamma$ and $\Sigma$ are smooth
$(N-1)$-dimensional submanifolds of $\partial \Omega$ with positive measures such that
$\Gamma \cap \Sigma = \emptyset$.
$\Delta_p$ is the $p$-Laplacian and $\frac{\partial}{\partial \nu}$ is the outer normal derivative. 
Here, $f$ is a subcritical perturbation of $|u|^{p^*-1}$. 

The sharpness of our result is estabished by the construction of suitable 
Palais-Smale sequences (corresponding to the critical level) which are not relatively compact.

Then we give the analogous condition to (\ref{seuil}) 
for a general system with critical exponents 
$$
\left\{
\begin{array}{lll}
-\Delta_{p}u &=& \lambda f(x,u) + u|u|^{\alpha-1} |v|^{\beta+1} \;\; \mbox{in} \;\; \Omega  \\ \\
-\Delta_{q}v &=& \mu g(x,v) + |u|^{\alpha+1} |v|^{\beta-1}v  \;\; \mbox{in} \;\; \Omega
\end {array}
\right. 
$$
together with Dirichlet or mixed boundary conditions, where $f$ and $g$ are subcritical perturbations of $|u|^{p^*-1}$ and 
$|v|^{q^*-1}$ respectively, $p^*=\frac{Np}{N-p}$ (resp. $q^*=\frac{Nq}{N-q}$) 
is the critical exponent of the Sobolev embedding $W^{1,p}(\Omega) \subset L^r(\Omega)$ 
(resp. $W^{1,q}(\Omega) \subset L^r(\Omega)$). Our approach provides a general condition based on the Nehari manifold, which can be 
extended to a large class of critical nonlinear problems. In this work, we confine ourselves to systems involving   
$(p,q)-$Laplacian operators and critical nonlinearities.  
The sharpness of our result is estabished, in the special case $p=q$, by the construction of suitable 
Palais-Smale sequences which are not relatively compact. The question of sharpness corresponding to the case $p  \ne q$ is still open.

For a more complete description of nonlinear elliptic systems, we 
refer the reader to the papers by De Figueiredo \cite{FIG1} and by De Figueiredo $\&$ Felmer \cite{FIG2} and the references therein.
\section{A general local compactness result}
For the reader's convenience, we start with the scalar case and to render the paper selfcontained we will recall or show some well-known facts.
\subsection{The scalar case}
Let $\Omega \subset \mathbb{R}^N$, $N \geq 3$, be a bounded domain with smooth boundary $\partial \Omega$. Let 
$f(x,u) \, : \; \Omega \times \mathbb{R} \to \mathbb{R}$ be a function which is measurable in $x$, continuous in $u$ and 
satisfying the growt condition at infinity
\begin{equation} \label{croissance}
|f(x,u)| = o(u^{p^*-1}) \;\; \mbox{as} \;\; u \to +\infty, \;\; \mbox{uniformly in} \;\; x. 
\end{equation}
This situation occurs, for example, in the special cases $f(x,u) = u$ or $f(x,u) = u^{q-1}$, $1<q<p^*$.

Consider the problem
\begin{equation} \label{gen_eq2}
\begin{array}{l}
-\Delta_p u = \lambda f(x,u) + |u|^{p^*-2}u  \;\; {\rm in } \;\; \Omega , \\ \\
u|_{\Gamma} = 0 \;\; \mbox{and} \;\; \frac{\partial{u}}{\partial \nu}|_{\Sigma}=0,
\end{array}
\end{equation}
where $\Omega$ is a bounded domain in $\mathbb{R}^N$, $N \geq 3$, 
with smooth boundary $\partial \Omega = \overline \Gamma \cup \overline \Sigma$,  
where $\Gamma$ and $\Sigma$ are smooth
$(N-1)$-dimensional submanifolds of $\partial \Omega$ with positive measures such that
$\Gamma \cap \Sigma = \emptyset$.
Problem (\ref{gen_eq2}) is posed in the framework of the Sobolev space 
$$
W^{1,p}_{\Gamma}(\Omega)= \{u \in W^{1,p}(\Omega) \; : \; u |_{\Gamma} = 0\}, \;\; 
$$
which is the closure of $C^1_0(\Omega \cap \Gamma,\mathbb{R})$
with respect to the norm of $W^{1,p}(\Omega)$. 
Notice that $meas (\Gamma) > 0$ implies that the Poincar\'e inequality is still available in 
$W^{1,p}_{\Gamma}(\Omega)$, so it can be endowed with the norm 
$$
||u||= || \nabla u ||_p
$$
and $(W^{1,p}_{\Gamma}(\Omega),|| \, . \, ||)$ is a reflexive and separable Banach space. The associated Euler-Lagrange functional 
is given by
$$
J_\lambda (u) := \frac{1}{p} ||\nabla u||_p^p - \frac{1}{p^*} ||u||_{p^*}^{p^*} - \lambda \int_\Omega F(x,u(x)) \, dx
$$
the corresponding Euler-Lagrange functional, where $F(x,u) := \int_0^u f(x,s) \, ds$. 

We recall here that the Nahari manifold associated to the functional $J_\lambda$ is given by:
$$
{\mathcal N}_{J_\lambda} = \{u \in W^{1,p}_\Gamma(\Omega) \setminus\{0\} \; : \; J_\lambda'(u)(u)=0\},
$$
and it is clear that ${\mathcal N}_{J_\lambda}$ contains all nontrivial critical points of $J_\lambda$. 
This manifold can be characterized more explicitely by the following
$$
{\mathcal N}_{J_\lambda} = \left\{tu, \; (t,u) \in  (\mathbb{R} \setminus\{0\}) \times (W^{1,p}_\Gamma(\Omega) \setminus\{0\} ) 
\; : \; \frac{d}{dt} J_\lambda(tu)=0\right\},
$$
where $t \mapsto J_\lambda(tu)$ is a function defined from $\mathbb{R}$ to itself, for every $u$ given in $W^{1,p}_\Gamma(\Omega) \setminus\{0\}$. 
We define the critical level associated to Problem (\ref{gen_eq2}) by:
\begin{equation} \label{cr_lev}
c^*(\lambda):=\inf_{w \in {\mathcal N}_{J_0}} J_0 (w) + 
\inf_{w \in {\mathcal N}_{J_\lambda} \cup \{0\}} J_\lambda (w).
\end{equation}
At this stage, we can state and show our first result
\begin{theorem} \label{general}
Let $\lambda \in \mathbb{R}$ and $(u_n)$ be a Palais-Smale sequence of $J_\lambda$ such that
\begin{equation} \label{elhamidi}
\lim_{n \to +\infty} J_{\lambda} (u_n) < c^*(\lambda).
\end{equation}
Then $(u_n)$ is relatively compact.
\end{theorem}
{\bf Proof. } 
Let $\lambda \in \mathbb{R}$ and  $(u_n)$ be a Palais-Smale sequence for $J_\lambda$ of level $c \in \mathbb{R}$ ((PS)$_c$ for short) 
satisfying the condition (\ref{elhamidi}). We claim that $(u_n)$ is bounded in $W^{1,p}_\Gamma(\Omega)$. Indeed, on has one hand
\begin{equation} \label{niveau}
\frac{1}{p} ||\nabla u_n||_p^p - \frac{1}{p^*} ||u_n||^{p^*}_{p^*} - \lambda \int_\Omega F(x,u_n) \, dx = c + o_n(1),
\end{equation}
and
\begin{equation} \label{derivee}
||\nabla u_n||_p^p - ||u_n||^{p^*}_{p^*} - \lambda \int_\Omega f(x,u_n)u_n \, dx = o_n(||\nabla u_n||_p).
\end{equation}
Then, 
$$
\left( \frac{1}{p} - \frac{1}{p^*} \right) ||u_n||^{p^*}_{p^*} + \frac{\lambda}{p} \int_\Omega f(x,u_n)u_n \, dx  - 
\lambda \int_\Omega F(x,u_n) \, dx = c + o_n(1) + o_n(||\nabla u_n||_p).
$$
Now, let $\varepsilon>0$, using the growth condition (\ref{croissance}), there exists
$c_1(\varepsilon)>0$ such that
$$
|f(x,u)| \leq \varepsilon |u|^{p^*-1} + c_1 \;\; \mbox{and} \;\; |F(x,u)| \leq \frac{\varepsilon}{p^*} |u|^{p^*} + c_1, \;\; 
\mbox{a.e.} \;\; x \in \Omega \;\; \mbox{and for every} \;\; u \in \mathbb{R}.
$$
Applying the H\"older and the Young inequalities to the last relations, it follows
\begin{equation} \label{combin}
||u_n||^{p^*}_{p^*}  \leq   \varepsilon ||\nabla u_n||_p + c_2(|\Omega|,\lambda,\varepsilon).
\end{equation}
Combining (\ref{combin}) and (\ref{niveau}), we deduce that $(u_n)$ is in fact bounded in $W^{1,p}_\Gamma(\Omega)$. So  
passing, if necessary to a subsequence, we can consider that 
\begin{eqnarray*}
u_n &\rightharpoonup& u \;\; \mbox{in} \;\; W^{1,p}_\Gamma(\Omega), \\
u_n &\to& u \;\; \mbox{a.e. in} \;\; \Omega.
\end{eqnarray*}
On the other hand, the growth condition (\ref{croissance}) implies also that, for almost every $x \in \Omega$, 
the functions $s \mapsto F(x,s)$ and 
$s \mapsto sf(x,s)$ satisfy the conditions of the Br\'ezis-Lieb Lemma (see Theorem 2 in \cite{BL}). 
Thus, we get the identities
\begin{eqnarray*}
\int_\Omega F(x,v_n) \, dx &=& \int_\Omega F(x,u_n) - \int_\Omega F(x,u) + o_n(1), \\
\int_\Omega f(x,v_n)v_n \, dx &=& \int_\Omega f(x,u_n)u_n - \int_\Omega f(x,u)u + o_n(1).
\end{eqnarray*}
Moreover, let $\varepsilon > 0$, there is $c_1(\varepsilon) > 0$ such that 
$$
\left | \int_\Omega f(x,v_n)v_n \, dx \right | \leq \varepsilon || v_n ||_{p^*}^{p^*} + c_1 ||v_n||_1.
$$
Let $C>0$ (which is independent of $n$ and $\varepsilon$), such that $|| v_n ||_{p^*}^{p^*} \leq C$.
Since $(v_n)$ converges strongly to $0$ in $L^1(\Omega)$, there is $n_0(\varepsilon) \in \mathbb{N}$ such that 
$||v_n||_1 \leq \varepsilon / c_1$, for every $n \geq n_0(\varepsilon)$, and consequently
$$
| \int_\Omega f(x,v_n)v_n \, dx | \leq \varepsilon (1 + C), \;\;\; \forall n \geq n_0(\varepsilon). 
$$
In the same way, rewriting $F(x,v_n) = \int_0^{v_n} f(x,s) \, ds$ and using the same arguments as above, we deduce that
\begin{eqnarray}
\int_\Omega F(x,v_n) \, dx&=& o_n(1) \label{Fvn}\\
\int_\Omega f(x,v_n)v_n \, dx &=& o_n(1) \label{fvnvn}.
\end{eqnarray}
Applying once again the Br\'ezis-Lieb Lemma, we conclude that $u \in {\mathcal N}_{J_\lambda} \cup \{0\}$ and  
\begin{eqnarray}
||v_n||^p - ||v_n||^{p^*}_{p^*} &=& o_n(1), \label{commune} \\
J_0 (v_n) := \frac{1}{p}||v_n||^p - \frac{1}{p^*}||v_n||^{p^*}_{p^*} &=& c - J_\lambda (u) + o_n(1). \label{Jzero}
\end{eqnarray}
A direct computation gives
$$
{\mathcal N}_{J_0} = \left\{ t_0(u) u \; : \; u \in W^{1,p}_\Gamma(\Omega) \setminus\{0\}
\right\},
$$ 
where
$$
t_0(u) := \left(\frac{||u||^p}{||u||_{p^*}^{p^*}} \right)^{\frac{1}{p^*-p}}.
$$
Now, let $b$ be the common limit of $||v_n||^p$ and $||v_n||_{p^*}^{p^*}$. Suppose that $b \ne 0$. On one hand we have
\begin{eqnarray*}
J_0(t_0(v_n)v_n) &=& 
\left( \frac{1}{p}-\frac{1}{p^*} \right) \left( \frac{||v_n||^p}{||v_n||_{p^*}^{p}}\right)^{\frac{p^*}{p^*-p}} \\
& \geq & \inf_{w \in {\mathcal N}_{J_0}} J_0 (w).
\end{eqnarray*}
Then 
$$
\lim_{n \to +\infty} J_0(t_0(v_n)v_n) = \frac{b}{N} \geq \inf_{w \in {\mathcal N}_{J_0}} J_0 (w). 
$$
On the other hand,  the identity (\ref{Jzero}) leads to
$$
\frac{b}{N} = c - J_\lambda (u).
$$
It follows then
\begin{eqnarray*}
c & \geq &  \inf_{w \in {\mathcal N}_{J_0}} J_0 (w) + J_\lambda (u) \\
  & \geq & \inf_{w \in {\mathcal N}_{J_0}} J_0 (w) +  \inf_{w \in {\mathcal N}_{J_\lambda} \cup \{0\}} J_\lambda (w), 
\end{eqnarray*}
which contradicts the condition (\ref{elhamidi}). This achives the proof.
\hfill $\Box$
\subsection{Sharpness of the critical level formula in the scalar case}
To show the sharpness of the critical level formula (\ref{elhamidi}), it suffices to carry out a Palais-Smale sequence for $J_\lambda$ 
of level $c^*(\lambda)$ which contains no convergent subsequence. \\
Consider, for a given $\varepsilon >0$, the extremal function 
$$
\Phi_{\varepsilon}(x)=C_{N} \varepsilon^{\frac{N-p}{p^{2}}}\Big( \varepsilon + |x|^{\frac{p}{p-1}} \Big)^{\frac{p-N}{p}} \;\; \mbox{with} \;\;
C_N:=\left( N \left( \frac{N-p}{p-1}\right) ^{p-1}\right)^{(N-p)/p^2}
$$
which attains the best constant $S$ of the Sobolev embedding
$$
D^{1,p}(\mathbb{R}^{N}) \hookrightarrow L^{p^{*}}(\mathbb{R}^{N}).
$$

Without loss of generality, we can consider that $0 \in \Sigma$.
Moreover, the set $\partial{\Omega}$ satisfies the following
property (see more details in Adimurthi, Pacella and Yadava \cite{APY}): \\
There exist $\delta >0$, an open neighborhood ${\mathcal V}$ of
$0$ and a diffeomorphism \linebreak $ \Psi \, : \, B_{\delta}(0)
\longrightarrow {\mathcal V}$ which has a jacobian determinant
equal to one at 0, with $\Psi (B_{\delta}^{{+}})={\mathcal V} \cap
\Omega$, where $B_{\delta}^{+}=B_{\delta}(0) \cap \{x \in
\mathbb{R}^N \, : \; x_N > 0 \}$. \\
Let $\varphi \in C_{0}^{\infty}(\mathbb{R}^N)$ such that $\varphi \equiv 1$ 
in a neighborhood of the origin. \\
We define the sequence defined by
\begin{equation} \label{sharpness}
\psi_{n}(x):=\varphi(x) \Phi_{1/n}(x), \;\;\; \mbox{for} \;\;\; n \in \mathbb{N}^*. 
\end{equation}
It is well known that the sequence $(\psi_n) \subset W^{1,p}_\Gamma(\Omega)$ is a Palais-Smale sequence for $J_0$ of level 
$\inf_{w \in {\mathcal N}_{J_0}} J_0 (w)$, which satisfies 
\begin{eqnarray*} 
\psi_n & \to & 0 \;\; \mbox{a.e. in} \;\; \Omega, \\
\nabla \psi_n & \to & 0 \;\; \mbox{a.e. in} \;\; \Omega, \\
||\psi_n||_{p^*}^{p^*} &\longrightarrow& \left[ N \, \inf_{w \in {\mathcal N}_{J_0}} J_0 (w) \right]^{p/N} := \ell \;\; \mbox{as} \;\; n \longrightarrow +\infty, \\
||\nabla \psi_n||_{p}^{p} &\longrightarrow& \left[ N \, \inf_{w \in {\mathcal N}_{J_0}} J_0 (w) \right]^{p/N} := \ell \;\; \mbox{as} \;\; n \longrightarrow +\infty.
\end{eqnarray*}
Now, let $(u_n)$ be a Palais-Smale sequence of $J_\lambda$ of level $\inf_{w \in {\mathcal N}_{J_\lambda} \cup \{0\}} J_\lambda (w)$. 
We will not go into further details concerning which subcritical terms $f(u)$ allow the existence of such sequences, but in the litterature,  
this occurs for various classes of subcritical terms. Applying Theorem \ref{general}, there exists a subsequence, still denoted by $(u_n)$, 
which converges to some $u \in W^{1,p}_\Gamma(\Omega)$. Then 
\begin{eqnarray*}
|| u_n + \psi_n ||_{p^*} & \leq & C, \\
u_n + \psi_n & \to & u \;\; \mbox{a.e. in} \;\; \Omega, \\
|| \nabla u_n + \nabla \psi_n ||_p & \leq & C, \\
\nabla u_n + \nabla \psi_n &\to & \nabla u \;\; \mbox{a.e. in} \;\; \Omega .
\end{eqnarray*}
where $C$ a positive constant independent of $n$. 
We apply the Br\'ezis-Lieb Lemma to the sequence $(u_n+\psi_n)$ and get 
$$
|| u_n + \psi_n ||_{p^*}^{p^*} = || (u_n-u) + \psi_n ||_{p^*}^{p^*} + ||u|| _{p^*}^{p^*} + \mbox{o}_n (1).
$$
Moreover, one has
$$ - || u_n-u ||_{p^*} + || \psi_n ||_{p^*} - \ell^{1/p^*} \leq || (u_n-u) + \psi_n ||_{p^*} - \ell^{1/p^*}  \leq  
|| u_n-u ||_{p^*} + || \psi_n ||_{p^*} - \ell^{1/p^*} $$
which implies that 
$$
|| (u_n-u) + \psi_n ||_{p^*} - \ell^{1/p^*} =  \mbox{o}_n (1).
$$
Therefore, we conclude that 
$$
|| u_n + \psi_n ||_{p^*}^{p^*} = ||u|| _{p^*}^{p^*} +  \ell + \mbox{o}_n (1).
$$
The same argumets applied to the sequence $(\nabla u_n + \nabla \psi_n)$ give
$$
|| \nabla u_n + \nabla \psi_n ||_{p}^{p} = ||\nabla u|| _{p}^{p} +  \ell + \mbox{o}_n (1).
$$
Finally, using the fact that 
\begin{eqnarray}
| \psi_n |^{p^*}  &\stackrel{*}{\rightharpoonup} & \ell \, \delta_0  \;\;\; \mbox{weakly}* \mbox{in} \;\; \mathcal{M}^+ (\Omega) \label{delta0}\\
| \nabla \psi_n |^p &\stackrel{*}{\rightharpoonup} & \ell \, \delta_0 \;\;\; \mbox{weakly}* \mbox{in} \;\; \mathcal{M}^+ (\Omega) \label{delta00}
\end{eqnarray}
where $\delta_0$ is the Dirac measure concentrated at the origin and $\mathcal{M}^+ (\Omega)$ is the space of positive finite measures \cite{WIL}),
we get that the sequence $(u_n + \psi_n)$ is a Palais-Smale sequence of $J_\lambda$ of level $c^*(\lambda)$.

We hence constructed a Palais-Smale sequence $(u_n + \psi_n)$ of $J_\lambda$ of level $c^*(\lambda)$ which can not be 
relatively compact in $W^{1,p}_\Gamma(\Omega)$. This justifies the sharpness of the critical level formula (\ref{elhamidi}).

\begin{remark} If we are interested by the homogeneous Dirichlet conditions, {\it i.e.} if $\Sigma = \emptyset$, 
the same arguments developed above are still valid, it suffices to assume that the origin $0 \in \Omega$ and 
consider $\varphi \in C_{0}^{\infty}(\Omega)$ such that $\varphi \equiv 1$ 
in a neighborhood of the origin.
\end{remark}

\subsection{The system case}
Now, consider the system 
\begin{equation} \label{syst}
\left\{
\begin{array}{lll}
-\Delta_{p}u &=& \lambda f(x,u) + u |u|^{\alpha-1} |v|^{\beta+1}, \\ \\
-\Delta_{q}v &=& \mu g(x,v) + |u|^{\alpha+1} |v|^{\beta-1}v,
\end{array}
\right.
\end{equation}
together with Dirichlet or mixed boundary conditions 
\begin{equation} \label{mix}
\left\{
\begin{array}{lll}
u|_{\Gamma_1} = 0 \;\;\; \mbox{and} \;\;\; \frac{\partial{u}}{\partial \nu}|_{\Sigma_1}=0, \\ \\
v|_{\Gamma_2} = 0 \;\;\; \mbox{and} \;\;\; \frac{\partial{v}}{\partial \nu}|_{\Sigma_2}=0,
\end{array}
\right.
\end{equation} 
where, $\Omega$ is a bounded domain in $\mathbb{R}^N$, $N \geq 3$, 
with smooth boundary $\partial \Omega = \overline \Gamma_i \cup \overline \Sigma_i$,  
where $\Gamma_i$ and $\Sigma_i$ are smooth
$(N-1)$-dimensional submanifolds of $\partial \Omega$ with positive measures such that
$\Gamma_i \cap \Sigma_i = \emptyset$, $i \in \{1,2\}$.
$\Delta_p$ is the $p$-Laplacian and $\frac{\partial}{\partial \nu}$ is the outer normal derivative. 
Also, it is clear that when $\Gamma_1=\Gamma_2=\partial \Omega$, one deals with homogeneous 
Dirichlet boundary conditions. We assume here that 
\begin{equation} \label{H1}
1 <  p < N, \;\; 1 < q < N,
\end{equation} 
and the critical condition
\begin{equation} \label{H2}
\frac{\alpha+1}{p^*} + \frac{\beta+1}{q^*} = 1.
\end{equation} 
Indeed, this condition represents the maximal growth such that the integrability of 
the product term $|u|^{\alpha+1} |v|^{\beta+1}$ (which will appear in the Euler-Lagrange functional) 
can be guaranteed by suitable H\"older estimates.

The functions $f$ and $g$ are two caratheodory functions which satisfy the growth conditions 
\begin{eqnarray} 
|f(x,u)| &=& o(u^{p^*-1}) \;\; \mbox{as} \;\; u \to +\infty, \;\; \mbox{uniformly in} \;\; x,  \label{croissancef}\\
|g(x,v)| &=& o(v^{q^*-1}) \;\; \mbox{as} \;\; v \to +\infty, \;\; \mbox{uniformly in} \;\; x.  \label{croissanceg}
\end{eqnarray}
Problem (\ref{syst}), together with (\ref{mix}), is posed in the framework of the Sobolev space 
$W=W^{1,p}_{\Gamma_1}(\Omega) \times W^{1,q}_{\Gamma_2}(\Omega)$, 
where 
$$
W^{1,p}_{\Gamma_1}(\Omega)= \{u \in W^{1,p}(\Omega) \; : \; u |_{\Gamma_1} = 0\}, \;\; 
W^{1,q}_{\Gamma_2}(\Omega)= \{u \in W^{1,q}(\Omega) \; : \; u |_{\Gamma_2} = 0\},
$$
which are respectively the closure of $C^1_0(\Omega \cap \Gamma_1,\mathbb{R})$
with respect to the norm of $W^{1,p}(\Omega)$ and $C^1_0(\Omega \cap \Gamma_2,\mathbb{R})$ with respect to the norm of 
$W^{1,q}(\Omega)$. 
Notice that $meas (\Gamma_i) > 0$, $i=1, \, 2$, imply that the Poincar\'e inequality is still available in 
$W^{1,p}_{\Gamma_1}(\Omega)$ and $W^{1,q}_{\Gamma_2}(\Omega)$, so $W$ can be endowed with the norm 
$$
||(u,v)||= || \nabla u ||_p + || \nabla v ||_q
$$
and $(W,|| \, . \, ||)$ is a reflexive and separable Banach space. The associated Euler-Lagrange functional 
$I_{\lambda, \mu} \in C^1(W,\mathbb{R})$ is given by 
$$
I_{\lambda, \mu} (u,v) = (\alpha+1)\left( \frac{P(u)}{p}  - \lambda \int_{\Omega} F(x,u) \right) + 
(\beta+1) \left( \frac{Q(v)}{q}  - \mu \int_{\Omega} G(x,v) \right) - R(u,v),
$$
where
$P(u)= ||\nabla u||_p^p,$
$Q(v)=||\nabla v||_q^q,$
$F(x,u) = \int_0^u f(x,s) \, ds,$
$G(x,v) = \int_0^v g(x,t) \, dt,$
and 
$R(u,v)=\int_\Omega |u|^{\alpha+1} |v|^{\beta+1} dx$. Notice that $R(u,v) \leq || u ||_{p^*}^{\alpha+1} || v ||_{q^*}^{\beta+1} < +\infty$.

Consider the Nehari manifold associated to Problem (\ref{syst}) given by
$$
\mathcal{N_{\lambda, \mu}} = \{ (u,v) \in W \setminus \{(0,0)\} \;  / \;
D_1 I_{\lambda, \mu} (u,v)(u) = D_2 I_{\lambda, \mu} (u,v)(v) = 0 \}, 
$$
where $D_1 I_{\lambda, \mu}$ and $D_2 I_{\lambda, \mu}$ are the derivative of $I_{\lambda, \mu}$ with respect to 
the first variable and the second variable respectively.

An interesting and useful characterization of  $\mathcal{N_{\lambda, \mu}}$ is the following
$$
\mathcal{N_{\lambda, \mu}} = \{(su,tv) \; / \; (s,u,t,v) \in \mathcal{Z^*} \; \mbox{and} \; 
\partial_{s} I_{\lambda, \mu} (su,tv) = \partial_{t} I_{\lambda, \mu} (su,tv) =0\},
$$
where 
$$
\mathcal{Z^*} = \{(s,u,t,v); \; (s,t) \in \mathbb{R}^2, \; (u,v) \in W^{1,p}_{\Gamma_1}(\Omega) \times W^{1,q}_{\Gamma_2}(\Omega), 
(su,tv) \ne (0,0)\} 
$$
and $I_{\lambda, \mu}$ is considered as a functional of four variables $(s,u,t,v)$ in 
$\mathcal{Z} := \mathbb{R} \times W^{1,p}_{\Gamma_1}(\Omega) \times \mathbb{R} \times W^{1,q}_{\Gamma_2}(\Omega)$. \\ 
%
\begin{definition} \label{palais_smale}
Let $\lambda$ and $\mu$ be two real parameters. 
A sequence $(u_n,v_n) \in W$ is a Palais-Smale sequence of the functional $I_{\lambda,\mu}$ if 
\begin{eqnarray} 
&\bullet& \mbox{there exists} \;\; c \in \mathbb{R} \;\; \mbox{such that} \;\; \lim_{n \to +\infty}I_{\lambda, \mu} (u_n,v_n) = c
\label{niv_c} \\
&\bullet& D I_{\lambda,\mu} (u_n,v_n) \;\; \mbox{converges strongly in the dual} \;\; W' \;\; \mbox{of} \;\; W \label{diff}
\end{eqnarray}
where $D I_{\lambda,\mu} (u_n,v_n)$ denotes the G\^ateaux derivative of $I_{\lambda,\mu}$.
\end{definition}
The last condition (\ref{diff}) implies that 
\begin{eqnarray} 
D_1 I_{\lambda,\mu} (u_n,v_n)(u_n) &=& \mbox{o} \, (|| u_n ||_{p^*}) 
\label{diff1} \\
D_2 I_{\lambda,\mu} (u_n,v_n)(v_n) &=& \mbox{o} \,(|| v_n ||_{q^*}).
\label{diff2}
\end{eqnarray}
where $D_1 I_{\lambda,\mu} (u_n,v_n)$ (resp. $D_2 I_{\lambda,\mu} (u_n,v_n)$) denotes 
the G\^ateaux derivative of $I_{\lambda,\mu}$ with respect to its first (resp. second) variable. \\

We introduce the critical level corresponding to Problem (\ref{syst}) by
\begin{equation} \label{cr_lev_sys}
c^*(\lambda,\mu):= \inf_{w \in {\mathcal N}_{0,0}} I_{0,0} (w) + 
\inf_{w \in \mathcal{N_{\lambda, \mu}} \cup \{(0,0)\}} I_{\lambda, \mu} (w).
\end{equation}
Then we have the following
\begin{theorem} \label{th1}
Let $\lambda$ and $\mu$ be two real parameters and $(u_n,v_n)$ be a Palais-Smale sequence of $I_{\lambda, \mu}$  such 
that 
\begin{equation} \label{elhamidi_sys}
c:=\lim_{n \to +\infty} I_{\lambda, \mu} (u_n,v_n) < c^*(\lambda,\mu).
\end{equation}
Then $(u_n,v_n)$ relatively compact.
\end{theorem}
{\bf Proof. } Let $\lambda$ and $\mu$ be two real parameters and $(u_n,v_n)$ be a Palais-Smale sequence of $I_{\lambda, \mu}$ 
satisfying the condition (\ref{elhamidi_sys}). We claim that $(u_n,v_n)$ is bounded in $W$. Indeed, on one hand conditions 
(\ref{niv_c}), (\ref{diff1}) and (\ref{diff2}) can be rewritten as the following
\begin{eqnarray} 
I_{\lambda, \mu} (u_n,v_n) &=& c + o_n (1) 
\label{niv} \\
P(u_n) - \lambda \int_{\Omega} f(x,u_n) u_n \, dx &=& R(u_n,v_n) + \mbox{o} \, (|| u_n ||_{p^*}) 
\label{der1} \\
Q(v_n) - \mu \int_{\Omega} f(x,v_n) v_n \, dx &=& R(u_n,v_n) + \mbox{o} \,(|| v_n ||_{q^*}).
\label{der2}
\end{eqnarray}
Using (\ref{H2}), one gets 
\begin{eqnarray} 
R(u_n,v_n) &=&  \frac{\alpha + 1}{p^*} \left( P(u_n) - \lambda \int_{\Omega} f(x,u_n) u_n \right) + \mbox{o} \, (|| u_n ||_{p^*}) 
\nonumber \\
           &+&  \frac{\beta + 1}{q^*} \left( Q(v_n) - \mu \int_{\Omega} g(x,v_n) v_n \right) + \mbox{o} \, (|| v_n ||_{q^*}).  
	   \label{Ruv}
\end{eqnarray}
Suppose that there is a subsequence, still denoted by $(u_n,v_n)$ in $W$ which is unbounded, {\it i.e.} 
$||\nabla u_n||_p + ||\nabla v_n||_q$ tends to $+\infty$ as $n$ goes to $+\infty$. \\
If 
$$
\lim_{n \to +\infty} ||\nabla u_n||_p = +\infty, 
$$
then using (\ref{croissancef}) one has 
\begin{eqnarray*} 
\int_{\Omega}  | f(x,u_n) u_n | &=& \mbox{o} \, (P(u_n)), \\
\int_{\Omega}  | F(x,u_n) | &=& \mbox{o} \, (P(u_n)),
\end{eqnarray*}
since (\ref{croissancef}) implies that for every $\varepsilon>0$, there exists 
$c_1(\varepsilon)>0$ such that
$$
|f(x,s)| \leq \varepsilon |s|^{p^*-1} + c_1 \;\; \mbox{and} \;\; |F(x,s)| \leq \frac{\varepsilon}{p^*} |s|^{p^*} + c_1, \;\; 
\mbox{a.e.} \;\; x \in \Omega, \;\; \forall \,  s \in \mathbb{R}.
$$
Similarly, if 
$$
\lim_{n \to +\infty} ||\nabla v_n||_q = +\infty, 
$$
then using (\ref{croissanceg}) it follows
\begin{eqnarray*} 
\int_{\Omega}  | g(x,v_n) v_n | &=& \text{o} \, (Q(v_n)), \\
\int_{\Omega}  | G(x,v_n) | &=& \mbox{o} \, (Q(v_n)).
\end{eqnarray*}
On one hand, suppose that
$$
\lim_{n \to +\infty} ||\nabla u_n||_p = \lim_{n \to +\infty} ||\nabla v_n||_q = +\infty. 
$$
Substituting (\ref{Ruv}) in (\ref{niv}), we obtain
\begin{eqnarray*} 
c + \mbox{o}_n (1) &=& (\alpha + 1) \left( \frac{1}{p} - \frac{1}{p^*} + \mbox{o} \, (P(u_n))^{\frac{p^*-p}{p}} \right) P(u_n) \\
                   &+& (\beta + 1) \left( \frac{1}{q} - \frac{1}{q^*} + \mbox{o} \, (Q(v_n))^{\frac{q^*-q}{q}} \right) Q(v_n) \;\;
		   \displaystyle{\longrightarrow_{n \to +\infty}} +\infty
\end{eqnarray*}
which can not hold true. On the other hand, suppose that  
$$
\lim_{n \to +\infty} ||\nabla u_n||_p = +\infty \;\; \mbox{and the sequence} \;\; ||\nabla v_n||_q \;\; \mbox{is bounded}, 
$$
then (\ref{der1}) implies that $R(u_n,v_n)$ is unbounded while (\ref{der2}) implies, on the contrary, that $R(u_n,v_n)$ is bounded. 
The case
$$
\lim_{n \to +\infty} ||\nabla v_n||_q = +\infty \;\; \mbox{and the sequence} \;\; ||\nabla u_n||_p \;\; \mbox{is bounded}, 
$$
leads to a contradiction with the same argument, which achieves the claim. \\
At this stage, we can assume, up to a subsequence, that 
\begin{eqnarray*}
u_n & \rightharpoonup & u \;\; \mbox{in} \;\; W^{1,p}_{\Gamma_1}(\Omega), \\ 
v_n & \rightharpoonup & v \;\; \mbox{in} \;\; W^{1,q}_{\Gamma_2}(\Omega), \\ 
u_n & \to & u \;\; \mbox{a.e. in} \;\; \Omega, \\
v_n & \to & v \;\; \mbox{a.e. in} \;\; \Omega.
\end{eqnarray*}
It is clear that 
$$
(u,v) \in \mathcal{N_{\lambda, \mu}} \cup \{(0,0)\}.
$$
Let us set 
$$
X_n = u_n-u \;\; \mbox{and} \;\; Y_n = v_n-v.
$$
Using again the growth conditions (\ref{croissancef}) and (\ref{croissanceg}), we show easily that the functions, which are defined on 
$\Omega \times \mathbb{R}$:
$(x,s) \mapsto sf(x,s)$, $(x,s) \mapsto sg(x,s)$, $(x,s) \mapsto F(x,s)$ and 
$(x,s) \mapsto G(x,s)$ satisfy the conditions of the Br\'ezis-Lieb lemma \cite{BL}. Then, we have the decompositions
\begin{eqnarray*}
\int_\Omega F(x,X_n)  &=& \int_\Omega F(x,u_n) - \int_\Omega F(x,u) + o_n(1), \\
\int_\Omega f(x,X_n)X_n  &=& \int_\Omega f(x,u_n)u_n - \int_\Omega f(x,u)u + o_n(1), \\
\int_\Omega G(x,Y_n)  &=& \int_\Omega G(x,v_n) - \int_\Omega G(x,v) + o_n(1), \\
\int_\Omega g(x,Y_n)Y_n  &=& \int_\Omega g(x,v_n)v_n - \int_\Omega g(x,v)v + o_n(1). 
\end{eqnarray*}
Moreover, let $\varepsilon > 0$, then there is $c_1(\varepsilon) > 0$ such that 
$$
\left | \int_\Omega f(x,X_n)X_n \, dx \right | \leq \varepsilon || X_n ||_{p^*}^{p^*} + c_1 ||X_n||_1.
$$
Let $C$ be a positive constant such that $|| X_n ||_{p^*}^{p^*} \leq C$.
Since $X_n$ converges to $0$ in $L^1(\Omega)$, there exists $n_0(\varepsilon) \in \mathbb{N}$ verifying 
$||X_n||_1 \leq \varepsilon / c_1$, for every $n \geq n_0(\varepsilon)$, thus
$$
\left| \int_\Omega f(x,X_n)X_n \, dx \right | \leq \varepsilon (1 + C), \;\;\; \forall n \geq n_0(\varepsilon). 
$$
In the same manner, writing $F(x,X_n) = \int_0^{X_n} f(x,s) \, ds$ and using the same arguments as above, we get
$$
\int_\Omega F(x,X_n)  = o_n(1) \;\; \mbox{and} \;\;
\int_\Omega f(x,X_n)X_n = o_n(1).
$$
Similarly, it follows that
$$
\int_\Omega G(x,Y_n) = o_n(1) \;\; \mbox{and} \;\;
\int_\Omega g(x,Y_n)Y_n = o_n(1).
$$
Applying a slightly modified version of the Br\'ezis-Lieb lemma \cite{MS}, one has 
$$
R(X_n,Y_n) = R(u_n,v_n) - R(u,v) + o_n(1).
$$
It follows that
\begin{eqnarray*}
P(X_n) - R(X_n,Y_n) = o_n(1), \\
Q(Y_n) - R(X_n,Y_n) = o_n(1), \\
I_{0,0} (X_n,Y_n) = c - I_{\lambda, \mu} (u,v) + o_n(1).
\end{eqnarray*}
Notice that the Nehari manifold associated to $I_{0,0}$ is given by
$$
\mathcal{N}_{0,0} = \left\{ (s_0(u,v) u, t_0(u,v)v); \; (u,v) \in W^{1,p}_{\Gamma_1}(\Omega) \times W^{1,q}_{\Gamma_2}(\Omega), \;
u \not\equiv 0, \, v \not\equiv 0 \right\},
$$ 
where 
$$
s_0(u,v) = \left[\frac{P(u) Q(v)^{\frac{r(\beta+1)}{q(\alpha+1)}}}{R(u,v)^{\frac{r}{\alpha+1}}}\right]^{\frac{1}{r-p}}, \;\;
t_0(u,v) = t(s_0(u,v)),
$$
and 
$$
r=\frac{(\alpha+1)q}{q-(\beta+1)} > p , \; \;\; 
t(s)=\left[\frac{R(u,v)}{Q(v)}\right]^{\frac{r}{q(\alpha+1)}} s^{\frac{r}{q}}.
$$
Let $\ell$ be the common limit of $P(X_n)$, $Q(Y_n)$ and $R(X_n,Y_n)$. We claim that $\ell = 0$. By contradiction, 
suppose that $\ell \ne 0$, then on one hand we get
\begin{eqnarray} \label{K}
I_{0,0} (s_0(X_n,Y_n) X_n, t_0(X_n,Y_n)Y_n) &=& (\alpha+1) \left( \frac{1}{p} - \frac{1}{r} \right) K(X_n,Y_n),\\
& \geq & \inf_{w \in {\mathcal N}_{0,0}} I_{0,0} (w),  \nonumber
\end{eqnarray}
where 
$$
K(X_n,Y_n) = \left[\frac{P(X_n)^{(\alpha+1)} Q(Y_n)^{(\beta+1)\frac{p}{q}}}{R(X_n,Y_n)^{p}}\right]^{\frac{r}{(\alpha+1)(r-p)}}.
$$
A direct computation shows that 
$$
\lim_{n \to +\infty} K(X_n,Y_n) = \ell,
$$
therefore
$$
\lim_{n \to +\infty} I_{0,0} (s_0(X_n,Y_n) X_n, t_0(X_n,Y_n)Y_n) = \ell (\alpha+1) \left( \frac{1}{p} - \frac{1}{r} \right).
$$
On the other hand, 
\begin{eqnarray*}
\lim_{n \to +\infty} I_{0,0} (X_n,Y_n) &=& \ell \left( \frac{\alpha+1}{p} + \frac{\beta+1}{q} - 1 \right) \\
&=& \ell (\alpha+1) \left( \frac{1}{p} - \frac{1}{r} \right).
\end{eqnarray*}
Hence, we obtain
$$
\ell (\alpha+1) \left( \frac{1}{p} - \frac{1}{r} \right) = c - I_{\lambda, \mu} (u,v),
$$
and consequently
\begin{eqnarray*}
c & \geq & \inf_{w \in {\mathcal N}_{0,0}} I_{0,0} (w) + I_{\lambda, \mu} (u,v) \\
  & \geq & \inf_{w \in {\mathcal N}_{0,0}} I_{0,0} (w) + 
\inf_{w \in \mathcal{N_{\lambda, \mu}} \cup \{(0,0)\}} I_{\lambda, \mu} (w).
\end{eqnarray*}
This leads to a contradiction with (\ref{elhamidi_sys}), then $\ell = 0$, which achieves the proof. 
\hfill $\Box$

\begin{remark}
1) In the scalar case, we obtain the analogous of Theorem \ref{th1}, the proof follows easily with the same arguments. 
We note here that if we consider the special case (\ref{bn_eq}), direct computations show that 
$$
\inf_{w \in {\mathcal N}_{0}} I_{0} (w) =  \frac{1}{N} S^{\frac{N}{2}}
\;\; \mbox{and} \;\; 
\inf_{w \in \mathcal{N_{\lambda}} \cup \{0\}} I_{\lambda} (w) = 0,
$$
which recovers the famous Br\'ezis-Nirenberg condition (\ref{seuil}). \\
2) It is clear that our condition (\ref{elhamidi}) or (\ref{elhamidi_sys}) can be extended to a large class of quasilinear or semilinear differential operators: 
Leray-Lions type operators, fourth-order operators. \\ 
3) Using the H\"older inequality in the denominator $R(u,v)$, we get 
\begin{equation} \label{important}
\inf_{(u,v) \in {\mathcal N}_{0,0}} I_{0,0} (u,v) \geq (\alpha+1) \left( \frac{1}{p} - \frac{1}{r} \right) 
\left[S_p S_q^{\frac{p(\beta+1)}{q(\alpha+1)}}\right]^{\frac{r}{r-p}},
\end{equation}
where $S_p$ (resp. $S_q$)
denotes the best Sobolev constant in the embedding $W^{1,p}_{\Gamma_1}(\Omega) \subset L^{p^*}(\Omega)$ 
(resp. $W^{1,q}_{\Gamma_2}(\Omega) \subset L^{q^*}(\Omega)$).
\end{remark}
We end this note by the following interesting relation arising in the special case $p=q$ and $\Gamma_1=\Gamma_2$.
\begin{prop} \label{propos}
Assume that $p=q>1$. Then, 
$$
\inf_{(u,v) \in {\mathcal N}_{0,0}} I_{0,0} (u,v) = \frac{p}{N-p} S_p^{\frac{N}{p}}.
$$
\end{prop}
{\bf Proof.} In the special case $p=q$, direct computations give 
$$
p^*= \alpha + \beta + 2 \;\;\; \mbox{and} \;\;\; (\alpha+1) \left( \frac{1}{p} - \frac{1}{r} \right) = \frac{p}{N-p}.
$$
Then, using (\ref{important}), we conclude that 
$$
\inf_{(u,v) \in {\mathcal N}_{0,0}} I_{0,0} (u,v) \geq \frac{p}{N-p} S_p^{\frac{N}{p}}.
$$
On the other hand, let $(u_n) \subset W^{1,p}_{\Gamma_1}(\Omega)$ be a minimizing sequence of $S_p$. Then using the identity (\ref{K}), we get 
\begin{eqnarray*}
\inf_{w \in {\mathcal N}_{0,0}} I_{0,0} (w) \leq I_{0,0} (s_0(u_n,u_n) u_n, t_0(u_n,u_n)u_n) & = & \frac{p}{N-p}
\left[\frac{|| \nabla u_n ||^p_p}{||u_n||^p_{p^*}}\right]^{\frac{rp^*}{(\alpha+1)(r-p)}} \\
& = & \frac{p}{N-p} \left[\frac{|| \nabla u_n ||^p_p}{||u_n||^p_{p^*}}\right]^{\frac{N}{p}}.
\end{eqnarray*}
It is clear that the last quantity goes to  $\displaystyle \frac{p}{N-p} S_p^{\frac{N}{p}}$
as $n +\infty$, which achieves the proof.
\hfill $\Box$

\begin{remark}
For the sharpness of the critical level (\ref{elhamidi_sys}), we define the sequence $\psi_{n}(x):=\varphi(x) \Phi_{1/n}(x)$ as in (\ref{sharpness}).
We consider then a Palais-Smale sequence $(u_n,v_n)$ for $J_{\lambda,\mu}$ of level  $\inf_{w \in \mathcal{N_{\lambda, \mu}} \cup \{(0,0)\}} I_{\lambda, \mu} (w)$. 
Following the same argumets developed in the scalar case and using Proposition \ref{propos}, we prove that the 
sequence $(u_n+\psi_{n},v_n+\psi_{n})$ is a Palais-Smale sequence for
$J_{\lambda,\mu}$ of level $c^*(\lambda,\mu)$ and which can not be relatively compact in $W$. 
This implies the sharpness of the critical level formula (\ref{elhamidi_sys}).
\end{remark}

$$
\begin{array}{ll}
\begin{array}{lllll}
\mbox{Khalid Adriouch} \\
\mbox{Laboratoire de Math. \& Applications,} \\
\mbox{Universit\'e de la Rochelle,} \\
\mbox{17042 La Rochelle, France.}\\
\mbox{e-mail: kadriouc@univ-lr.fr}
\end{array} &
\begin{array}{lllll}
\mbox{Abdallah El Hamidi} \\
\mbox{Laboratoire de Math. \& Applications,} \\
\mbox{Universit\'e de la Rochelle,} \\
\mbox{17042 La Rochelle, France. }\\
\mbox{e-mail: aelhamid@univ-lr.fr}
\end{array}
\end{array}
$$

\end{document}